
\documentstyle{amsppt}

\magnification=1200

\def\O{{\Cal O}}
\def\C{{\Bbb C}}
\def\lra{\longrightarrow}
\def\ra{\rightarrow}
\def\x{\otimes}
\def\End{{\Cal E\text{\rm nd}}}
\def\Hom{{\Cal H\text{\rm om}}}
\def\hom{{\text{\rm Hom}}}
\def\ker{{\text{\rm ker}}}
\def\Diff{{\Cal D\text{\rm iff\,}}}
\def\sign{{\text{\rm sign}}}
\def\sym{{\text{\rm Sym}}}
\def\tot{{\text{\rm Tot}}}
\def\gr{{\text{\rm gr}}}
\def\Spec{{\text{\bf Spec}}}
\def\g{\frak g}

\newcount\headcount\newcount\subheadcount
\newcount\subsubheadcount\newcount\itemcount
\headcount=0\subheadcount=0
\subsubheadcount=0\itemcount=0

\def\headno{
\global\subheadcount=0\global\subsubheadcount=0
\global\advance\headcount by1 \number\headcount. }
\def\subheadno{
\global\subsubheadcount=0\global\advance\subheadcount by1
\number\headcount.\number\subheadcount. }
\def\itemno{\global\advance\itemcount by1 \number\itemcount }
\def\namething#1#2{\csname newbox\endcsname#1\setbox#1=\hbox{#2}\copy#1}

\def\makeyoungdiagram#1{\vbox
{\offinterlineskip\ialign{\rlap
{\hbox to12pt{\hfill##\hfill}}\rlap
{\vrule height9.5ptdepth2.5ptwidth.2pt}\rlap
{\smash{\vrule height-2.5ptdepth2.7ptwidth12.2pt}}\vrule
height9.5ptdepth-9.3ptwidth12pt\rlap
{\vrule height9.5ptdepth2.5ptwidth.2pt}&&\rlap
{\hbox to12pt{\hfill##\hfill}}\rlap
{\vrule height9.5ptdepth2.5ptwidth.2pt}\rlap
{\smash{\vrule height-2.5ptdepth2.7ptwidth12.2pt}}\vrule
height9.5ptdepth-9.3ptwidth12pt\rlap
{\vrule height9.5ptdepth2.5ptwidth.2pt}\cr
#1}}}

\topmatter

\title
Resolution for Sheaf of Differential Operators on Smooth Free Geometric
Quotient of Linear Action of Algebraic Group
\endtitle
\rightheadtext{Differential Operators on Quotient}

\author
Gwoho Liu
\endauthor

\date
July 24, 1996
\enddate

\affil
Department of Mathematics\\ University of California\\
Riverside, CA 92521, USA
\endaffil

\address
Department of Mathematics, University of California,
Riverside, CA 92521, USA
\endaddress

\email
st20c\@math.ucr.edu
\endemail

\endtopmatter

\document

In [R], Z.~Ran gave a canonical construction for the universal
deformation of a simple vector bundle using the Jacobi complex
of an appropriate differential graded Lie algebra.
Independently, H.~Esnault and E.~Viehweg made a similar construction.
Using these tools, we obtain a resolution for the sheaf of differential
operators on smooth geometric quotients of free linear actions of
algebraic groups.  Unlike previous applications of these methods,
we are able to obtain global results.

The terms of our resolution involve symmetric and alternating powers
of vector bundles easily constructed geometrically from the algebraic
group and the vector space on which it acts.
Our resolution is particularly simple for projective space.
For $\Bbb P(V)$, we obtain the following sequence---
$$
0 \lra \O \lra S^r(V^*(1)) \lra \Diff^r_{\Bbb P(V)}/\O \lra 0,
$$
in which the first map is symmetrization of the map in the Euler sequence.
We are able to conclude that the sheaf of differential operators on
$\Bbb P^n$, with $n>1$, has no non-scalar global endomorphisms.

\medskip

This paper forms part of the author's UCR Ph.~D.~dissertation.

We work over the complex numbers.

\head \headno Schur Functor of a Complex of Vector Spaces
\endhead

We define the Schur functor of a complex of vector spaces.
We provide sign computations because such seems to be rare,
though we will not need to use them in much generality.
Signs were actually computed from scratch, so the signs contained here
can be considered to be independent experimental data.
Mostly this parallels the description of Schur functors given in [F].

\subhead \subheadno How to Get a Symmetrizer from a Partition
\endsubhead

We describe how one constructs an element of the group ring, $\C[\sym(n)]$
,from a partition of $n$.  This is exactly as in [F].

\definition {Definition}
A {\it partition} of an integer, $n$, is a non-increasing sequence of
positive integers, summing to $n$.
\enddefinition

We fix a partition of $n$, $l = (l_i)$, for the remainder
of this section.

\medskip

Given our partition we draw a picture of boxes, called the Young diagram.
We put $l_i$ boxes in the $i^{\text{\rm th}}$ row,
and we label them with integers from one to $n$.
For example, for the partition $(4,2,1)$, we draw the following picture---
\medskip
\makeyoungdiagram {
	1&2&3&4\cr
	5&6\cr
	7\cr
}

\medskip

Let $S < \sym(n)$ be the subgroup of the symmetric group consisting of
all permutations preserving the rows of the Young diagram,
and let $L$ be the subgroup preserving columns.
We define the Young symmetrizer to be the following
element of $\C[\sym(n)]$---
$$
\Sigma(l) = \left(\sum_{\sigma\in S}\sigma\right)
\left(\sum_{\sigma\in L}\sign(\sigma)\sigma\right).
$$

For example, for the diagram,
\makeyoungdiagram {
	1&2\cr
	3\cr
},
we have $\Sigma(2,1) = (1+(12))(1-(13)) = (1+(12)-(13)-(132))$.

\medskip

Notice that if we have a $\sym(n)$ action on a $\C$--vector space,
we get a $\C[\sym(n)]$ action.

\subhead \subheadno Schur Functor for Complexes of Vector Spaces
\endsubhead

We first define the Schur functor of a graded vector space.
Let $l$ be a partition of $n$.
Let $A = \oplus A_i$ be a graded vector space.
We define an action of a permutation, $\sigma\in\sym(n)$,
on $A^{\x n}$ by the following formula---
$$
a_1\x\ldots\x a_n \mapsto
\sign(\sigma,\deg a_1,\ldots,\deg a_n)
a_{\sigma 1}\x\ldots\x a_{\sigma n},
$$
in which
$$
\sign(\sigma,d_1,\ldots,d_n) = 
\prod_{i<j\atop\sigma j<\sigma i}
(-1)^{d_i d_j}.
$$
If $o$ of the $d_i$ are odd, we obtain a homomorphism from $\sym(n)$ to
$\sym(o)$ by forgetting the even things.
By viewing $\sign(\sigma,d_1,\ldots,d_n)$ as the sign of the
image of $\sigma$ under this homomorphism, we see that
$\sign(\cdot,d_1,\ldots,d_n)$ is a homomorphism, thus we do indeed
obtain action.

\definition {Definition}
Given $A$ and $l$ as above,
the {\it Schur module}, $S^l A$, is the image of the
Young symmetrizer, $\Sigma(l): A^{\x n} \ra A^{\x n}$.
\enddefinition

Since $S^l A$ is naturally a direct summand of $A^{\x n}$,
we see that $S^l(\cdot)$ is a functor, which we call the Schur functor.

\medskip

We now extend our definition to case of a complex.
If $A$ is a complex with differential $d$,
we make $A^{\x n}$ into a complex with differential being
$$
a_1\x\ldots\x a_n\mapsto
\sum_i \Big(\prod_{j<i}(-1)^{\deg a_j}\Big)
a_1\x\ldots\x da_i\x\ldots\x a_n.
$$

\medskip

\proclaim{Proposition \itemno}
$\Sigma$ and $d$ commute.
\endproclaim

\demo{proof}
We prove this with direct computation.
Let $\Sigma = \sum_{\sigma\in\sym(n)}c(\sigma)\sigma$,
with $c(\sigma)\in\C$.
The image of $a_1\x\ldots\x a_n$ under $d\circ\Sigma$ is
$$
\multline
\sum_i\Big(\prod_{j<i}(-1)^{\deg a_j}\Big)
\sum_{\sigma\in\sym(n)}
c(\sigma)
\\
\sign(\sigma,\deg a_1,\ldots,\deg a_i+1,\ldots,\deg a_n)
a_{\sigma1}\x\ldots\x da_i\x\ldots\x a_{\sigma n}.
\endmultline
$$
We change the order of the sums and sum over $\sigma i$ in place of $i$ and
obtain the following---
$$
\multline
\sum_{\sigma\in\sym(n)}
\sum_i\Big(\prod_{j<\sigma i}(-1)^{\deg a_j}\Big)
c(\sigma)
\\
\sign(\sigma,\deg a_1,\ldots,\deg a_{\sigma i}+1,\ldots,\deg a_n)
a_{\sigma1}\x\ldots\x da_{\sigma i}\x\ldots\x a_{\sigma n}.
\endmultline
$$
The image of $a_1\x\ldots\x a_n$ under $\Sigma\circ d$ is
$$
\multline
\sum_{\sigma\in\sym(n)}
c(\sigma)
\sign(\sigma,\deg a_1,\ldots,\deg a_n)
\sum_i
\\
\Big(\prod_{j<i}(-1)^{\deg a_{\sigma j}}\Big)
a_{\sigma1}\x\ldots\x da_{\sigma i}\x\ldots\x a_{\sigma n}.
\endmultline
$$
Thus we need to compare, for fixed $\sigma$ and $i$,
$$
\sign(\sigma,\deg a_1,\ldots,\deg a_n)
\Big(\prod_{j<i}(-1)^{\deg a_{\sigma j}}\Big)
$$
with
$$
\sign(\sigma,\deg a_1,\ldots,\deg a_{\sigma i}+1,\ldots,\deg a_n)
\Big(\prod_{j<\sigma i}(-1)^{\deg a_j}\Big).
$$
We see that they are equal because
$$
\multline
\sign(\sigma,\deg a_1,\ldots,\deg a_n)
\sign(\sigma,\deg a_1,\ldots,\deg a_{\sigma i}+1,\ldots,\deg a_n)
\\
= \sign(\sigma^{-1},\deg a_1,\ldots,\deg a_n)
\sign(\sigma^{-1},\deg a_1,\ldots,\deg a_{\sigma i}+1,\ldots,\deg a_n)
\\
= \Big(\prod_{\sigma i<\sigma j\atop j<i}(-1)^{\deg a_{\sigma j}}\Big)
\Big(\prod_{\sigma j<\sigma i\atop i<j}(-1)^{\deg a_{\sigma j}}\Big)
= \Big(\prod_{j<\sigma i}(-1)^{\deg a_j}\Big)
\Big(\prod_{j<i}(-1)^{\deg a_{\sigma j}}\Big)
.
\endmultline
$$
\qed
\enddemo

Commutativity of $d$ and $\Sigma$ tells us that $S^l(A)$ is
a direct summand of the complex $A^{\x n}$ as a complex.
In particular, $S^l(\cdot)$ takes quasi--isomorphic complexes
to quasi--isomorphic complexes.

The dual of a partition is the partition corresponding the the
transpose of its Young diagram.
For example, the dual of $(3,1)$ is $(2,1,1)$.
If $a$ and $b$ are dual partitions, we write
$\wedge^b = S^a$.

\head \headno Definitions of DGLA and Jacobi Complex
\endhead

Given a complex, $L$, with a map $[,]:\wedge^2L \ra L$, we obtain a map
$L\x[,]:L\x\wedge^2L \ra L\x L$.
Symmetrizing we obtain a map $\wedge^3L \ra \wedge^2L$.
Similarly we obtain maps $\wedge^iL \ra \wedge^{i-1}L$.
Given another complex, $M$, with a map
$\langle,\rangle:L\x M\ra M$, we obtain a map
$L\x\langle,\rangle:L\x L\x M\ra L\x M$.
Symmetrizing we obtain a map $\wedge^2L\x M\ra L\x M$.

\definition {Definition}
A {\it differential graded Lie algebra} (DGLA) is a complex,
$L$, together with a map $[,]:\wedge^2L \ra L$ such that the
composition $\wedge^3L \ra \wedge^2L \ra L$ is zero.
A {\it module} of $L$ is a complex, $M$, with a map
$\langle,\rangle:L\x M\ra M$ such that the following square commutes---
$$
\CD
\wedge^2L\x M @>{[,]\x M}>> L\x M \\
@V{L\x\langle,\rangle}VV @V{\langle,\rangle}VV \\
L\x M @>{\langle,\rangle}>> M \\
\endCD
$$
\enddefinition

\definition {Definition}
Let $L$ be a DGLA.
{\it Jacobi complex}, $J^r(L)$, is the following complex---
$$
\CD
\wedge^rL @>>> \ldots @>>> \wedge^3L @>>> \wedge^2L @>>> L,
\endCD
$$
with $\wedge^iL$ in degree $-i$ for $1\le i\le r$ and zeros elsewhere.
The {\it stupid} filtration of $J^r(L)$ is the filtration,
$F^iJ^r(L) = J^i(L)$ for $i\le r$.
\enddefinition

Notice that if $L$ is a DGLA and $M$ an $L$--module, then
the complex $L {\buildrel0\over\lra} M$, with the zero map, naturally
forms a DGLA as follows---
$$
\CD
S^2M \\
@A0AA \\
L\x M @>{\langle,\rangle}>> M \\
@A0AA @A0AA \\
\wedge^2L @>{[,]}>> L. \\
\endCD
$$
$J^r(L\ra M)$ splits up as a direct sum of complexes,
with a summand for each degree in $M$.

\definition {Definition}
Let $M$ be a $L$--module.
The other {\it Jacobi complex}, $J^r(L,M)$, is the summand of
$J^r(L{\buildrel0\over\lra}M)$ having degree one in $M$.
\enddefinition

\definition {Definition}
If $V$ is a finite dimensional vector space, we call
$V$ a {\it coalgebra} if $V^*$ is an algebra.
In this case, if $W$ is another a finite dimensional vector space, we call
$W$ a $V$ {\it module} if $W$ is a $V^*$ module.
\enddefinition

The natural map $J^r(L) \ra S^2J^r(L)$,
$$
\matrix
\ldots & \lra &
\wedge^3L\x L\oplus S^2\wedge^2L & \lra & \wedge^2L\x L & \lra & \wedge^2L
\\
&& \uparrow && \uparrow && \uparrow
\\
\ldots & \lra &
\wedge^4L & \lra & \wedge^3L & \lra & \wedge^2L & \lra & L,
\endmatrix
$$
gives a map
$h^0(J^r(L)) \ra S^2h^0(J^r(L))$ and the map
$J^r(L,M) \ra J^r(L)\x J^r(L,M)$ gives a map
$h^0(J^r(L,M)) \ra h^0(J^r(L))\x h^0(J^r(L,M))$.
We call these maps of cohomology groups comultiplication.
This comultiplication turns
$h^0(J^r(L))\oplus\C$ into a coalgebra and $h^0(J^r(L,M))$ into a module.

\definition {Definition}
Let $V$ and $W$ be vector bundles on a smooth algebraic variety, $X$.
Let $X\times X$ be the product with projections $p_i$,
and diagonal $\Delta$.
The {\it Sheaf of differential operators of order $r$ from $V$ to $W$},
$\Diff^r(V,W)$, is $\Hom(p_{1*}(p_2^*(V)\x\O_{X\times X}/I_\Delta^{r+1}),W)$.
We write $\Diff^r$ for $\Diff^r(\O,\O)$.
\enddefinition

The ring structure on $\O_{X\times X}$ makes
$\Diff^r$ into a sheaf of $\O_X$--coalgebras, and
$\Diff^r(V,W)$ into a $\Diff^r$--module.

\head \headno The Result
\endhead

Let $X$ be an algebraic variety and $G$ be an algebraic group
bundle over $X$.  Let $\g$ be the corresponding Lie algebra bundle.
Let $V$ be a vector bundle and
let $A:G\times_X V\ra V$ be a linear action.
Let $s\in\Gamma(V)$ be a section.  From
this section, we obtain a map $G\ra V$.
Let $a:\g\ra V$ be the (relative) tangent of the map,
$G\ra V, g\mapsto g(s)$,
and let $b:\g\x V\ra V$ be map one obtains from the (relative) tangent
of the map $G\ra GL(V)$.

We wish to understand something of the map $G\times V\ra V$
in terms of $a$ and $b$.

\proclaim{Lemma \itemno}
Using notation as above,
let $\{1\}_1 = \Spec(\O_X\oplus\g^*)$ be the first infinitesimal
neighborhood of $1$ in $G$
and $\{s\}_r$ be the $r^{\text{\rm th}}$ infinitesimal neighborhood
of $s$ in $V$.  The map
$$
f:\bigoplus_{i=0}^{r+1}S^iV^* \lra
\bigoplus_{i=0}^rS^iV^*\bigoplus\bigoplus_{i=0}^r\g^*\x S^iV^*,
$$
corresponding to $\{1\}_1 \times \{s\}_r \ra \{s\}_{r+1}$,
is given by
$$
\Pi_{i=1}^dv_i\mapsto \Pi_{i=1}^dv_i\bigoplus
a^*(v_1)\x\Pi_{i=2}^dv_i+b^*(v_1)\x\Pi_{i=2}^dv_i.
\eqno{(*)}
$$
\endproclaim

\demo{proof}
We think of $\{1\}_1$ as an affine algebraic group scheme over $X$,
with structure sheaf $\O_X\oplus\g^*$.
We identify $\{s\}_r$ with $\{0\}_r = \Spec(\oplus_{i=0}^{r+1}S^iV^*)$
by translation by $s$, which we denote $t_s: V\ra V$.
The question is local on $X$, so we let $\{e_i\}$ be a local basis of $V$,
and let $\{x_i\}$ be the dual basis.
The map in question is $t_{-s}\circ A\circ(G\times t_s): G\times V\ra V$.
Since $(*)$ preserves multiplication and $V^*$ generates $\sym V^*$,
we only need to compute the map from $V^*$.
Let $v\in \Gamma(V^*)$.
Let $c:V\x V^*\ra\g^*$ be the dual of the tangent of the map
$G\ra GL(V)$.
$t_{-s}^\#$, the map of sheaves of $\O_X$ algebras corresponding to
$t_{-s}$, takes $v$ to $v - v(s)$.
$A^\#$ takes $v-v(s)$ to $\sum_i c(v\x e_i)\x x_i + v - v(s)$.
Finally $(G\times t_s)^\#$ takes this to
$\sum_i c(v\x e_i)\x x_i + x_i(s)c(v\x e_i) + v + v(s) - v(s)$.
\qed
\enddemo

Let $U\subset V$ be an open union of $G$ orbits.
Assume that a geometric quotient $\pi: U\ra Y$ exists and
that $Y$ is smooth over $X$.
Assume that the section, $s$, is contained in $U$.
Let $p:Y\ra X$.
Let $I_s$ be the ideal sheaf of the image of $s$ in $Y$.

\proclaim{Proposition \namething\propone\itemno}
We use notation as above.
Give $L = (\g\ra V)$ the natural DGLA structure.
Then $h^0(J^r(L))^*\oplus\O = p_*(\O_Y/I_s^{r+1})$.
If the map, $\g\ra V$, is injective, then $J^r(L)$ has no other cohomology.
Also the multiplication on $p_*(\O_Y/I_s^{r+1})$ coming from
multiplication on $\O_Y$ coincides with the natural multiplication
on $h^0(J^r(L))^*$.
\endproclaim

\demo{proof} From
the map $\pi:U\ra Y$ we get a map of sheaves of $\O_X$ algebras---
$$
\pi^\#:p_{1*}(O_Y/I_s^{r+1}) \ra
\oplus_{i=0}^rS^iV^*.
$$
The action, $A:G\times V\ra V$, and
projection, $p_2$, give us two $\O_X$ algebra maps---
$$
A^\#,p_2^\#:
\oplus_{i=0}^rS^iV^* \ra
\oplus_{i=0}^rS^iV^* \bigoplus
\oplus_{i=1}^rS^{i-1}V^*\x\g^*.
$$
$p_2^\#$ is identity on the right summand of the target and zero on
the left summand.
Lemma identifies $A^\#$.
Thus $A^\# - p_2^\#$ maps to the right summand only.
Notice that $A^\# - p_2^\#$ coincides with the dual of the final map
of $\tot(J^r(L))$.

We make a sequence, which we will show to be left exact,
out of $A^\# - p_2^\#$ and $\pi^\#$---
$$
0 \lra
p_{1*}(O_Y/I_s^{r+1}) \lra
\oplus_{i=0}^rS^iV^* \lra
\oplus_{i=1}^rS^{i-1}V^*\x\g^*.
$$
$(A^\# - p_2^\#)\circ\pi^\#$ is zero because $A^\#$ and $p_2^\#$ agree on
$G$ invariant functions.
$\pi^\#$ is injective because $\pi$ is surjective.
Thus $p_{1*}(O_Y/I_s^{r+1})$ maps injectively
to $\ker(A^\#-p_2^\#) = h^0(J^r(L))^* \oplus \O$.
This map is filtered because both range and domain have the filtration
induced from the filtration on $\oplus_{i=0}^rS^iV^*$.

Because $h^1(L) = T_{Y/X}$ and $L$ has no higher cohomology, we see that
$h^i(\wedge^iL) = S^iT_{Y/X}$ and $\wedge^iL$ has no higher cohomology.
This tells us that
$$
\gr(h^0(J^r(L))) = \oplus_{i=1}^rS^iT_{Y/X}.
$$
Therefore the map from $p_{1*}(O_Y/I_s^{r+1})$ to $h^0(J^r(L))^*\oplus\O$
is isomorphism on all graded parts.  Thus our sequence is left exact.

In the case that $\g\ra L$ is injective, $\wedge^iL$ has no cohomology
other than $h^i$ therefore $J^r(L)$ has no cohomology other than $h^0$.

Since $\pi^\#$ preserves multiplication, the multiplication on
$p_{1*}(O_Y/I_s^{r+1})$ is the one inherited from
$\oplus_{i=0}^rS^iV^*$.  This latter multiplication corresponds to the
comultiplication on $h^0(J^r(V[-1]))$, in which $V[-1]$ is given the
trivial DGLA structure.
Since $V[-1]\ra L$ is a map of DGLA,
$h^0(J^r(L))$ has comultiplication inherited from $h^0(J^r(V[-1]))$.
This proves that the multiplications coincide.
\qed
\enddemo

Let us fix some notation.
Let $G$ be an algebraic group and $V$ be a vector space.
Let $A:G\times V\ra V$ be a free linear action.
Let $\g = T_1(G)$.
Let $U\subset V$ be an open union of orbits.
Assume that a geometric quotient $\pi: U\ra X$ exists and that $X$ is smooth.

Let $G$ act on $G\times U$ by conjugation on the
first factor and $A$ on the second factor.
Let $\bar G$ be the quotient of $G\times U$.
By doing group arithmetic on the first factor and
mapping to $X$ by projection to the second factor,
we make $\bar G$ into an algebraic group bundle.
Let $\bar\g$ be the Lie algebra bundle.

Let $\bar V$ be the quotient of $V\times U$ by the diagonal action of $G$.
Notice that $\bar V$ is a $\bar G$--representation over $X$.
The diagonal map $U \ra V\times U$ gives a section, $\Delta$, of $\bar V$.
Let $\bar A:\bar G\times_X \bar V\ra \bar V$ be the action.
Let $\bar U$ be the quotient of $U\times U$.
Notice that the quotient of $\bar U$ by action of $\bar G$ is
$X\times X$ and the section, $\Delta$, maps to diagonal.

{}From the section $\Delta:X \ra \bar U$, we obtain a map $\bar G\ra\bar U$.
Let $a:\bar\g\ra\bar V$ be the (relative) tangent of this map.
Let $b:\bar\g\x\bar V\ra\bar V$ be the map coming from the
(relative) tangent of the map $\bar G\ra GL(\bar V)$.

\proclaim{Corollary \itemno}
Give $L = (\bar\g\ra\bar V)$ the natural DGLA structure.
Then $h^0(J^r(L))\oplus\O = \Diff^r$ and $J^r(L)$ has no other
cohomology.  Also the coalgebra structures coincide.
\endproclaim

Let $W$ be another vector bundle and
$B:G\times W\ra W$ be a linear action.
As before, we make a $\bar G$ representation, $\bar W$, with action
$\bar B:\bar G\times_X\bar W\ra\bar W$.

\proclaim{Corollary \itemno}
Give $L = (\bar\g\ra\bar V)$ the natural DGLA structure and
give $\bar W$ the natural $L$--module structure.
Then $h^0(J^r(L,\bar W)) = \Diff^r(\bar W^*,\O_X)$ and
$J^r(L,\bar W)$ has no other cohomology.
Also the module structures coincide.
\endproclaim

\demo{proof}
We apply Proposition $\copy\propone$ to the $\bar G$ action on
$\bar U\oplus\bar W^*$, with the section $(\Delta,0)$.
Notice that the quotient is the geometric vector bundle
$p_2^*(\bar W^*)$ over $X\times X$.

We proceed according to the following general formula:
given the affine coordinate ring of a geometric vector
bundle, we obtain the associated module by extracting the degree one
component and dualizing;  module structure may be obtained from the
multiplication of the of degree zero component by the degree one
component of the affine coordinate ring.

Proposition tells us that the affine coordinate ring of
the infinitesimal neighborhood of $(\Delta,0)$ in $p_2^*(\bar W^*)$
is given by $h^0(J^{r+1}(L\ra \bar W^*))^* \oplus \O$.
The degree one component is $h^0(J^r(L,\bar W^*))^*$ and
the degree zero component is
$h^0(J^r(L,\O))^* = \O_{X\times X}/I_\Delta^{r+1}$.
We dualize the former with respect to the latter.
Since
$$
\Hom_{\O_{X\times X}}(p_2^*(\cdot),\O_{X\times X}) =
p_2^*\Hom_{\O_X}(\cdot,\O_X),
$$
we obtain
$$
h^0(J^r(L,\bar W))^* =
p_{1*}(p_2^*(\bar W^*)\x\O_{X\times X}/I_\Delta^{r+1}).
$$
Applying $\Hom(\cdot,\O_X)$ to both sides, we obtain
$$
h^0(J^r(L,\bar W)) = \Diff^r(\bar W^*,\O).
$$

We conclude that the module structures coincide by noticing that, in general,
the comultiplication of the degree zero by the degree one components of
$h^0(J^r(L\ra M))$ coincides with the comultiplication
$h^0(J^r(L,M)) \ra h^0(J^r(L))\x h^0(J^r(L,M))$.
\qed
\enddemo

\head \headno Application to Projective Space
\endhead

Let $V$ be a $n+1$ dimensional vector space,
and $Q$ be a one dimensional vector space.
We view $\Bbb P(V)$ as the quotient of $\hom(V,Q)\setminus\{0\}$ by $GL(Q)$.
Let $L = (\O \ra V^*(1))$.  We notice that
$$
J^r(L) =
\pmatrix
S^r(V^*(1)) \\
\uparrow \\
S^{r-1}(V^*(1)) &{\buildrel\sim\over\lra} \\
&&\ldots \\
&&& {\buildrel\sim\over\lra}& S^2(V^*(1)) \\
&&&& \uparrow \\
&&&& V^*(1) &{\buildrel\sim\over\lra}& V^*(1) \\
&&&&&& \uparrow \\
&&&&&& \O
\endpmatrix
$$
We simplify this and obtain the following short exact sequence---
$$
0 \lra \O \lra S^r(V^*(1)) \lra \Diff^r/\O \lra 0.
$$

Using this we identify
$$
H^0\left(\End(\Diff^r/\O)\right) =
\Bbb H^0\left(\End^\cdot(\O\ra S^r(V^*(1)))\right).
$$
The spectral sequence for hypercohomology
tells us that for $n>1$, this is one dimensional.
Thus $\Diff^r_{\Bbb P(V)}/\O$ has no non scalar global endomorphisms.

However, this is not the case for $\Diff^r(\O(k),\O(k))$ for $k\neq 0$.
Notice that $J^r(L,\O(k)) =$
$$
\matrix
S^r(V^*(1))(k) \\
\uparrow \\
S^{r-1}(V^*(1))(k) &{\buildrel k+r-1\over\lra} \\
&& \ldots \\
&&& {\buildrel k+2\over\lra} & S^2(V^*(1))(k) \\
&&&& \uparrow \\
&&&& V^*(k+1) &{\buildrel k+1\over\lra}& V^*(k+1) \\
&&&&&& \uparrow \\
&&&&&& \O(k) & {\buildrel k\over\lra}& \O(k) \\
\endmatrix
$$
We see that
$$
\Diff^r(\O(k),\O(k)) =
\cases
S^k(V^*(1)) \oplus S^r(V^*(1))/S^k(V^*(1))
& \text{\rm for~} 0\le k \le r \text{\rm ~and} \\
S^r(V^*(1)) & \text{\rm otherwise.}
\endcases
$$

\head References
\endhead
\tolerance=10000

\roster

\item"[R]"  Z.~Ran: `On the Local Geometry of Moduli Spaces of Locally Free
Sheaves', Moduli of Vector Bundles, Marcel Dekker, New York, 1996.

\item"[F]"  W.~Fulton and J.~Harris: `Representation Theory: A First Course'.
Springer-Verlag, 1991.

\endroster

\enddocument